\UseRawInputEncoding
\documentclass[11pt,onecolumn]{amsart}
\newtheorem{theorem}{Theorem}[section]
\usepackage[left=1.5cm,top=1.5cm,bottom=1.5cm,right=1.5cm]{geometry}
\usepackage{latexsym}
\usepackage{amssymb,amsmath,amsfonts}
\usepackage[pagewise]{lineno}
\usepackage{multirow}
\usepackage{graphicx}
\usepackage[colorlinks = true]{hyperref}
\theoremstyle{definition}

\newtheorem{example}[theorem]{Example}

\newcommand{\N}{\mathbb N}

\numberwithin{equation}{section}
 \global\parskip 6pt \oddsidemargin
0.1cm \evensidemargin 0.4cm \headheight 0.1cm \topmargin -2.0cm
\begin{document}
\pagenumbering{arabic}
\title[Solving delay differential equations via Sumudu transform]
{Solving delay differential equations via Sumudu transform}
\author{ M.O. Aibinu$^{1*}$, S. C. Thakur$^2$, S. Moyo$^3$}
\address{$^{1}$ Institute for Systems Science \& KZN e-Skills CoLab, Durban University of Technology, Durban 4000, South Africa}
\address{$^{1}$  DSI-NRF Centre of Excellence in Mathematical and Statistical Sciences (CoE-MaSS), South Africa}
\address{$^{2}$ KZN e-Skills CoLab, Durban University of  Technology, Durban 4000, South Africa}
\address{$^{3}$Institute for Systems Science \& Office of the DVC Research, Innovation \& Engagement, Milena Court, Durban University of Technology, Durban 4000, South Africa}
\email{${^1*}$moaibinu@yahoo.com / mathewa@dut.ac.za} 
\email{$^2$thakur@dut.ac.za}
\email{$^3$dvcrie@dut.ac.za}

\keywords{Delay equations, Sumudu transform, Pantograph-type equations.\\
{\rm 2010} {\it Mathematics Subject Classification}: 34K13; 34K05; 06F30; 46B20;\\
{\bf Citation}: M. O. Aibinu, S.C. Thakur, S. Moyo, Solving delay differential equations via Sumudu transform, International Journal of Nonlinear Analysis and Applications, DOI: 10.22075/ijnaa.2021.22682.2402 (In Press).}

\begin{abstract}
A technique which is known as Sumudu Transform Method (STM) is studied for the construction of solutions of a most general form of delay differential equations of pantograph type. This is a pioneer study on using the STM to construct the solutions of delay differential equations of pantograph type with variable coefficients. We obtain the exact and approximate solutions of nonlinear problems with multiproportional delays and variable coefficients. The strength of STM is illustrated in reducing the complex computational work as compared to the well-known methods. This paper shows how to succinctly identify the Lagrange multipliers for nonlinear delay differential equations with variable coefficients, using the STM. The potency and suitability of the STM are exhibited by giving expository examples. The solutions of nonlinear Volterra integro-differential equations of pantograph type are also obtained.
\end{abstract}

\maketitle
\section{Introduction}
The important roles which Delay Differential Equations (DDEs) play in explaining many occurrences in physical and biological systems are enormous. Delays are introduced in models to enhance their vitality and to allow an accurate narration of actual problems which involve the use of DDEs. Modeling of the spread of infectious diseases, population dynamics or electrodynamics problems are governed by DDEs. DDEs contain terms whose value depends on the solution at prior times. Time delays can be a constant, time-dependent, or state-dependent. Several events in nature are not instantaneous in their occurrences and can therefore be modeled by DDEs. The study of DDEs has attracted much attention. There are reports on the conditions for the solutions and qualitative behaviour of solutions of considered differential equations \cite{Mohammed1, Moyo1, Graef, Moyo3, Moyo2}. Several approximation and numerical methods have been developed for solving DDEs. Perturbation methods top the list in the approximate methods which are being used for solving differential equations, algebraic equations,  integrodifferential equations and difference equations (See e.g, \cite{Nayfeh}). The impediment to these methods is the indispensability of a small parameter which might be indistinguishable from the original physical parameter of the given equations. Several research efforts were made to overcome this drawback. Diverse methods were introduced which include variational iteration method \cite{He1, He2, He3, He4, He5}, Laplace transform \cite{Wu1, Wu2}.

\par In this work, an innovative approach of variational iteration method via the Sumudu transform is studied for the construction of solutions of a most general form of DDEs of pantograph type. We consider the DDEs of pantograph type with variable coefficients of the form,
 \begin{eqnarray}\label{sumud17}
&&g\left(x, y(q_0x), y'(q_1x), y''(q_2x)\right)=0, x_0\in [x_0,T]\nonumber\\
&&y(x_0)=a, y'(x_0)=b,
 \end{eqnarray}
 where $q_i \in (0,1),$ for $i=1,2,$ introduced in 1971 \cite{Ockendon1}. The efficacy of the Sumudu Transform Method (STM) is shown in obtaining the exact and approximate solutions of nonlinear problems with multi-proportional delays and variable coefficients. The flexibility, consistency and effectiveness of the STM for solving DDEs with variable coefficients and nonlinear Volterra integro-differential equations are demonstrated. Expository examples are presented to show that there is no need for linearization, perturbations or large computational works. 
 
\section{Variational iteration method and Sumudu transform} Variational iteration method (VIM) was instituted by He (See e.g,  \cite{He4} and references there in). The flexibility, consistency, and effectiveness of the VIM  as compared to other well-known methods have been demostrated (See e.g, \cite{Wu1, Wu2} and references there in). Moreover, efforts have been made on the development of the VIM with a blend of Laplace transform (See e.g, Wu \cite{Wu1, Wu2}), and Sumudu transform for obtaing the solutions of DDEs with constant coefficients (See e.g, \cite{Vilu1}). Applications of Sumudu transform for obtaining the solutions of several forms of differential equations are also found in \cite{Hussein, Golmankhaneh, Alomari, Nisar}. The Sumudu transform is a simple modified form of Laplace transfom. A corresponding property of the Sumudu transform can be derived from every property which is true for Laplace transform and vice versa (See e.g, \cite{Watugala1, Belgacem1, Belgacem2}). The technique of Sumudu transform is a simple, effective and universal way for identifying the Lagrange multiplier. Sumudu transform is applicable provided the given function $g(x)$ satisfies the following Dirichlet conditions:
\begin{itemize}
	\item [(i)] it is single valued function which may have a finite number of finite isolated discontinuies for $x>0.$
	\item [(ii)] it remains less than $be^{-a_0x}$ as $x$ approaches $\infty,$ where $b$ is a positive constant and $a_0$ is a real positive number.
\end{itemize} 
\subsection{Presentation of Sumudu transform}
The concept of Sumudu Transform (ST) was proposed in 1993 for solving differential equations and control engineering problems by Watugala \cite{Watugala1}. ST is an integral transform. The ST of a function $g(x)$ is denoted as $G(u)$ and is defined for all real numbers $x\geq 0$ by
 \begin{equation}\label{sumud18}
G(u)=S\left[g(x)\right]=\int^{\infty}_{0}g(ux)e^{-x}dx.
 \end{equation}
 For the integer order derivatives, the ST is given as 
 \begin{equation}\label{sumud19}
S\left[\frac{dg(x)}{dx}\right]=\frac{1}{u}\left[G(u)-g(0)\right].
 \end{equation}
The ST for the $n$-order derivative is
 \begin{equation}\label{sumud20}
S\left[\frac{d^ng(x)}{dx^n}\right]=\frac{1}{u^n}\left[G(u)-\displaystyle\sum_{k=0}^{n-1}u^k\frac{d^kg(x)}{dx^k}|_{x=0}\right].
 \end{equation}
The ST is linear and it preserves linear functions. Also the ST is credited for its units preserving property, which makes it applicable for solving problems without the need to resort to a new frequency domain (See e.g, Watugala \cite{Watugala1}, Belgacem et al. \cite{Belgacem1}). The fundamental properties of ST are elucidated in Belgacem \& Karaballi \cite{Belgacem2}. The highlight of the ST is presented by considering the universal nonlinear equation,
 \begin{equation}\label{sumud21}
\frac{d^ny(x)}{dx^n}+R\left[y(x)\right]+N\left[y(x)\right]=f(x),
 \end{equation}
 subject to the initial conditions
 \begin{equation}\label{sumud22}
y^{(k)}(0)=a_k,
 \end{equation}
 where $y^{(k)}(0)=\frac{d^ky(0)}{dx^k},$ $k=0, 1, ..., n-1,$ $R$ is a linear operator, $N$ is a nonlinear operator, $f(x)$ is a given continuous function and the highest order derivative is $\frac{d^ny(x)}{dx^n}.$
 \par Let $Y(u)=S[y(x)],$ taking the ST of (\ref{sumud21}) transforms its linear part into an algebraic equation of the form
 \begin{equation}\label{sumud23}
\frac{1}{u^n}Y(u)-\displaystyle\sum_{k=0}^{n-1}\frac{1}{u^{n-k}}y^{(k)}(0)=S\left[f(x)-R\left[y\right]-N\left[y\right]\right].
 \end{equation}
 Thus, the corresponding iteration procedure is given by
  \begin{equation}\label{sumud24}
Y_{n+1}(u)=Y_n(u)+{\varphi}(u)\left(\frac{1}{u^n}Y_n(u)-\displaystyle\sum_{k=0}^{n-1}\frac{1}{u^{n-k}}y^{(k)}(0)-S\left[f(x)-R\left[y\right]-N\left[y\right]\right]\right),
 \end{equation}
 where ${\varphi}(u)$ is the Lagrange multiplier. Take the classical variation operator on both sides of (\ref{sumud24}) and consider $S\left[R\left[y\right]+N\left[y\right]\right]$ as restricted terms. It is obtained that 
  \begin{equation}\label{sumud25}
\delta Y_{n+1}(u)=\delta Y_{n}(u)+{\varphi}(u)\frac{1}{u^n}\delta Y_n(u),
 \end{equation}
 which gives 
 \begin{equation}\label{sumud26}
{\varphi}(u)=-u^n.
 \end{equation}
 Substituting (\ref{sumud26}) into (\ref{sumud24}) and taking the inverse-Sumudu transform $S^{-1}$ of (\ref{sumud24}) yields the explicit iterative procedure, 
  \begin{eqnarray}\label{sumud27}
y_{n+1}(x)&=&y_n(x)+S^{-1}\left[-u^n\left(\frac{1}{u^n}Y(u)-\displaystyle\sum_{k=0}^{n-1}\frac{1}{u^{n-k}}y^{(k)}(0)-S\left[f(x)-R\left[y\right]-N\left[y\right]\right]\right)\right],\nonumber\\
&=&y_1(x)+S^{-1}\left[u^n\left(S\left[f(x)-R\left[y(x)\right]-N\left[y(x)\right]\right]\right)\right],
 \end{eqnarray}
 where 
   \begin{eqnarray}\label{sumud28}
y_1(x)&=&S^{-1}\left[\displaystyle\sum_{k=0}^{n-1}u^{k}y^{(k)}(0)\right]\nonumber\\
&=&y(0)+y'(0)x+...+\frac{y^{n-1}(0)x^{n-1}}{(n-1)!}.
 \end{eqnarray}
\subsection{Variable coefficient nonlinear equation}
Suppose the universal nonlinear equation (\ref{sumud21}) is endowed with variable coefficients such that it reads
 \begin{equation}\label{sumud29}
 \frac{d^ny(x)}{dx^n}+\alpha R_1[y(x)] + \beta (x)R_2[y(x)]+N\left[y(x)\right]=f(x),
 \end{equation}
 where $\alpha$ is a constant and $\beta (x)$ is a variable coefficient, $R_1$ and $R_2$ are linear operators and other terms remain as defined in (\ref{sumud21}). Taking the ST of (\ref{sumud29}) generates the iteration procedure
   \begin{eqnarray}\label{sumud30}
Y_{n+1}(u)&=&Y_n(u)+{\varphi}(u)(\frac{1}{u^n}Y_n(u)-\displaystyle\sum_{k=0}^{n-1}\frac{1}{u^{n-k}}y^{(k)}(0)\nonumber\\
&&-S\left[f(x)-\alpha R_1[y] - \beta (x)R_2[y]-N\left[y\right]\right]).
 \end{eqnarray}
Then the Lagrange multiplier ${\varphi}(u)$ is derived with $S\left[\beta (x)R_2[y]+N\left[y\right]\right]$ as the restricted terms. The rest of the computation follows the same process.
\section{Solution of pantograph type equations}
Pantograph Type Equations (PTEs) are functional differential equations with proportional delays. An example of a generalized multi-pantograph type equation is given by (\ref{sumud17}). The name pantograph was derived in 1971 from the work of Ockendon and Tayler \cite{Ockendon1}. PTEs are pragmatic in modeling of several problems which occur in diverse fields such as biology, electrodynamics,  economy, cell growth, quantum mechanics, astrophysics, number theory, probability theory, nonlinear dynamical systems and other industrial applications. Interested readers are referred to the following references for more notes on the applications of pantograph type equations \cite{Ockendon1, Ajello1, Buhmann}. PTEs are kind of DDEs and many authors have made research efforts for obtaining their solutions, both numerically and analytically \cite{Bashi1, Zhan, Sezer1, Sezer2, Saadatmandi1}. 
 \subsection{Delay differential equations}
\begin{figure}
\includegraphics[width=12.0cm ,height=10.0cm]{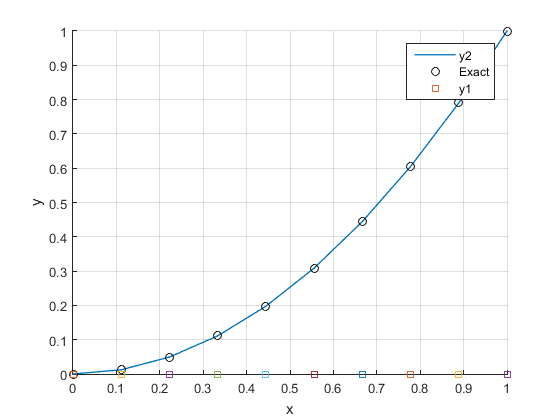}
 \caption{Graph of $y(x)=x^2.$}
 \label{sumud31}
\end{figure}

\begin{figure}
\includegraphics[width=12.0cm ,height=10.0cm]{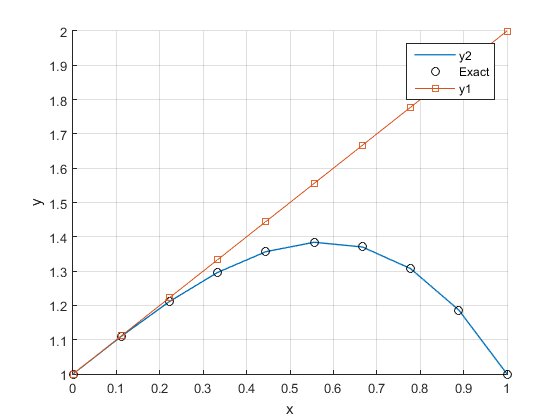}
 \caption{Graph of $y(x)=1+x-x^3.$}
 \label{sumud32}
\end{figure}

 Three examples are given to illustrate the iteration of STM for solving DDEs of multi-pantograph type with variable coefficients. After the first iteration, each subsequent iteration is used to obtain each succeeding term while the Higher Other Terms (HOT) are truncated. The Matlab is used to compare the graphs of the exact solutions with the solutions obtained by using the STM.
 \begin{example}\label{sumud2}
 Consider the nonlinear second-order pantograph equation
 \begin{eqnarray}\label{sumud12}
 y''(x)&-&8y\left(\frac{x}{2}\right)+xy'(x)=2,\nonumber\\
 y(0)&=&0,\\
 y'(0)&=&0.\nonumber
 \end{eqnarray}
\par {\bf Solution:} Taking the ST of (\ref{sumud12}) gives
\begin{equation}\label{sumud41}
\frac{Y(u)}{u^2}-\frac{y(0)}{u^2}-\frac{y'(0)}{u}=S\left[8y\left(\frac{x}{2}\right)-xy'(x)+2\right].
\end{equation}
Since $y(0)=0$ and $y'(0)=0,$ equation (\ref{sumud41}) gives
\begin{equation}
\frac{Y(u)}{u^2}=S\left[8y\left(\frac{x}{2}\right)-xy'(x)+2\right].
\end{equation}
Thus, the variational iteration formula is given by
\begin{eqnarray}\label{sumud3}
Y_{n+1}(u)&=&Y_n(u)+{\varphi}(u)\left(\frac{Y_n(u)}{u^2}-S\left[8y_n\left(\frac{x}{2}\right)-xy'_n(x)+2\right]\right), n\in \N.
\end{eqnarray}
Take the classical variation operator on both sides of (\ref{sumud3}) and consider the terms $S\left[8y\left(\frac{x}{2}\right)\right]$ and $S\left[xy'(x)\right]$ as the restricted variations. The Lagrange multiplier is obtained as
\begin{equation}
{\varphi}(u)=-u^2.
\end{equation}
Taking the inverse-Sumudu transform $S^{-1},$ of (\ref{sumud3}) gives the explicit iteration formula as
\begin{eqnarray}\label{sumud4}
y_{n+1}(x)&=&y_n(x)+S^{-1}\left[-u^2\left(\frac{Y_n(u)}{u^2}-S\left[8y_n\left(\frac{x}{2}\right)-xy'_n(x)+2\right]\right)\right]\nonumber\\
&=&y_1(x)+S^{-1}\left[u^2\left(S\left[8y_n\left(\frac{x}{2}\right)-xy'_n(x)+2\right]\right)\right],
\end{eqnarray}
where $y_1(x)$ is an initial approximation of (\ref{sumud4}) and which is obtained as $y_1(x)=0.$ Then,
\begin{eqnarray}\label{sumud10}
y_{2}(x)&=&S^{-1}\left[u^2\left(S\left[8y_1\left(\frac{x}{2}\right)-xy'_1(x)+2\right]\right)\right]\nonumber\\
&=&S^{-1}\left[u^2\left(S\left[2\right]\right)\right]\nonumber\\
&=&S^{-1}\left[2u^2\right]\nonumber\\
&=&x^2.
\end{eqnarray}
\par The exact solution is $y(x)=x^2.$ Figure \ref{sumud31} displays the exact and approximate solutions for Example \ref{sumud2}.
\end{example}
 \begin{example}\label{sumud5}
 Consider the nonlinear pantograph equation of second order
 \begin{eqnarray}\label{sumud13}
 y''(x)&-&\frac{8}{3}y'\left(\frac{x}{2}\right)y(x)-8x^2y\left(\frac{x}{2}\right)=-\frac{4}{3}-\frac{22}{3}x-7x^2-\frac{5}{3}x^3,~t\in [0,1],\nonumber\\
 y(0)&=&1,\\
 y'(0)&=&1.\nonumber
 \end{eqnarray}
\par {\bf Solution:} By taking the ST of (\ref{sumud13}), we obtain
\begin{equation}
\frac{Y(u)}{u^2}-\frac{y(0)}{u^2}-\frac{y'(0)}{u}=S\left[\frac{8}{3}y'\left(\frac{x}{2}\right)y(x)+8x^2y\left(\frac{x}{2}\right)-\frac{4}{3}-\frac{22}{3}x-7x^2-\frac{5}{3}x^3\right].
\end{equation}
Since $y(0)=1$ and $y'(0)=1,$ it is obtained that
\begin{equation}
\frac{Y(u)}{u^2}-\frac{1}{u^2}-\frac{1}{u}=S\left[\frac{8}{3}y'\left(\frac{x}{2}\right)y(x)+8x^2y\left(\frac{x}{2}\right)-\frac{4}{3}-\frac{22}{3}x-7x^2-\frac{5}{3}x^3\right].
\end{equation}
Thus for $ n\in \N,$ the variational iteration formula is given by
\begin{eqnarray}\label{sumud6}
Y_{n+1}(u)&=&Y_n(u)\nonumber\\
&+&{\varphi}(u)\left(\frac{Y_n(u)}{u^2}-\frac{1}{u^2}-\frac{1}{u}-S\left[\frac{8}{3}y'_n\left(\frac{x}{2}\right)y_n(x)+8x^2y_n\left(\frac{x}{2}\right)-\frac{4}{3}-\frac{22}{3}x-7x^2-\frac{5}{3}x^3\right]\right).
\end{eqnarray}
The classical variation operator on both sides of (\ref{sumud6}) is taken and the terms $S\left[\frac{8}{3}y'_n\left(\frac{x}{2}\right)y_n(x)\right]$ and $\left[8x^2y_n\left(\frac{x}{2}\right)\right]$ are being considered as the restricted variations. The Lagrange multiplier is obtained as
\begin{equation}
{\varphi}(u)=-u^2.
\end{equation}
The inverse-Sumudu transform $S^{-1},$ of (\ref{sumud6}) is taken to get the explicit iteration formula
\begin{eqnarray}\label{sumud7}
y_{n+1}(x)&=&y_n(x)\nonumber\\
&+&S^{-1}\left[-u^2\left(\frac{Y_n(u)}{u^2}-\frac{1}{u^2}-\frac{1}{u}-S\left[\frac{8}{3}y'_n\left(\frac{x}{2}\right)y_n(x)+8x^2y_n\left(\frac{x}{2}\right)-\frac{4}{3}-\frac{22}{3}x-7x^2-\frac{5}{3}x^3\right]\right)\right]\nonumber\\
&=&y_1(x)+S^{-1}\left[u^2\left(S\left[\frac{8}{3}y'_n\left(\frac{x}{2}\right)y_n(x)+8x^2y_n\left(\frac{x}{2}\right)-\frac{4}{3}-\frac{22}{3}x-7x^2-\frac{5}{3}x^3\right]\right)\right],
\end{eqnarray}
where the initial approximation of (\ref{sumud7}) is given by $y_1(x)=S^{-1}[1+u]=1+x$ and
\begin{eqnarray}\label{sumud8}
y_{2}(x)&=&y_1(x)\nonumber\\
&+&S^{-1}\left[u^2\left(S\left[\frac{8}{3}y'_1\left(\frac{x}{2}\right)y_1(x)+8x^2y_1\left(\frac{x}{2}\right)-\frac{4}{3}-\frac{22}{3}x-7x^2-\frac{5}{3}x^3\right]\right)\right].
\end{eqnarray}
Notice that $y_1\left(\frac{x}{2}\right)=1+\frac{x}{2}$ and $y'_1\left(\frac{x}{2}\right)=\frac{1}{2}.$ Recall that the task is to get the next term only. Therefore
\begin{eqnarray}\label{sumud9}
y_{2}(x)&=&1+x+S^{-1}\left[u^2\left(S\left[\frac{8}{3}\left(\frac{1}{2}\right)(1+x)-\frac{4}{3}-\frac{22}{3}x \right]\right)\right]~(\mbox{Truncate HOT})\nonumber\\
&=&1+x+S^{-1}\left[u^2\left(S\left[-6x\right]\right)\right]\nonumber\\
&=&1+x+S^{-1}\left[u^2(-6u)\right]\nonumber\\
&=&1+x+S^{-1}\left[-6u^3\right]\nonumber\\
&=&1+x-x^3.
\end{eqnarray}
\par The exact solution is known to be $y(x)=1+x-x^3.$ Figure \ref{sumud32} displays the exact and approximate solutions for Example \ref{sumud5}.
\end{example}

 \begin{example}\label{sumud11}
Consider the first order nonlinear DDE
 \begin{eqnarray}\label{sumud14}
 y'(x)-2xy^4\left(\frac{x}{2}\right)=0, ~y(0)=1.
 \end{eqnarray}
\par {\bf Solution:} The ST of (\ref{sumud14}) gives
 \begin{equation}
\frac{Y(u)}{u}-\frac{y(0)}{u}=S\left[2xy^4\left(\frac{x}{2}\right)\right].
\end{equation}
Since $y(0)=1,$ it is realized that 
\begin{equation}
\frac{Y(u)}{u}-\frac{1}{u}=S\left[2xy^4\left(\frac{x}{2}\right)\right].
\end{equation}
Then the variational iteration formula is given by
\begin{eqnarray}\label{sumud15}
Y_{n+1}(u)&=&Y_n(u)+{\varphi}(u)\left(\frac{Y_n(u)}{u}-\frac{1}{u}-S\left[2xy^4_n\left(\frac{x}{2}\right)\right]\right).
\end{eqnarray}
Consider $S\left[2xy^4_n\left(\frac{x}{2}\right)\right]$ as the restricted variation and take the classical variation operator on both sides of (\ref{sumud15}). This gives the Lagrange multiplier to be
\begin{equation}
{\varphi}(u)=-u.
\end{equation}
By taking the inverse-Sumudu transform $S^{-1}$ of (\ref{sumud15}), it gives the explicit iteration formula as
\begin{eqnarray}\label{sumud16}
y_{n+1}(x)&=&y_n(x)+S^{-1}\left[-u\left(\frac{Y(u)}{u}-\frac{1}{u}-S\left[2xy^4_n\left(\frac{x}{2}\right)\right]\right)\right]\nonumber\\
&=&y_1(x)+S^{-1}\left[u\left(S\left[2xy^4_n\left(\frac{x}{2}\right)\right]\right)\right],
\end{eqnarray}
with the initial approximation which is given by
$$y_1(x)=S^{-1}\left[-u(\frac{-1}{u})\right]=S^{-1}[1]=1.$$
\begin{eqnarray*}
y_{2}(x)&=&y_1(x)+S^{-1}\left[u\left(S\left[2xy^4_1\left(\frac{x}{2}\right)\right]\right)\right]\\
&=&1+S^{-1}\left[u\left(S\left[2x\right]\right)\right]\\
&=&1+S^{-1}\left[u(2u)\right]\\
&=&1+S^{-1}\left[2u^2\right]\\
&=&1+x^2.
\end{eqnarray*}
Next is to evaluate $y_3(x)$ to obtain the next term only.
\begin{eqnarray*}
y_{3}(x)&=&y_1(x)+S^{-1}\left[u\left(S\left[2xy^4_2\left(\frac{x}{2}\right)\right]\right)\right]\\
&=&1+S^{-1}\left[u\left(S\left[2x\left(1+\frac{x^2}{4}\right)^4\right]\right)\right]\\
&=&1+S^{-1}\left[u\left(S\left[2x\left(1+x^2\right)\right]\right)\right]~(\mbox{Truncate HOT})\\
&=&1+S^{-1}\left[u\left(S\left[2x+2x^3\right]\right)\right]\\
&=&1+S^{-1}\left[u\left(S\left[2x+2x^3\right]\right)\right]\\
&=&1+S^{-1}\left[u\left(2u+12u^3\right)\right]\\
&=&1+S^{-1}\left[2u^2+12u^4\right]\\
&=&1+x^2+\frac{x^4}{2}.
\end{eqnarray*}
Next is to evaluate $y_4(x)$ to obtain the next term only.
\begin{eqnarray*}
y_{4}(x)&=&y_1(x)+S^{-1}\left[u\left(S\left[2xy^4_3\left(\frac{x}{2}\right)\right]\right)\right]\\
&=&1+S^{-1}\left[u\left(S\left[2x\left(1+\frac{x^2}{4}+\frac{x^4}{32}\right)^4\right]\right)\right]\\
&=&1+S^{-1}\left[u\left(S\left[2x\left(1+x^2+\frac{x^4}{2}\right)\right]\right)\right]~(\mbox{Truncate HOT})\\
&=&1+S^{-1}\left[u\left(S\left[2x+2x^3+x^5\right]\right)\right]\\
&=&1+S^{-1}\left[u\left(2u+12u^3+120u^5\right)\right]\\
&=&1+S^{-1}\left[2u^2+12u^4+120u^6\right]\\
&=&1+x^2+\frac{x^4}{2}+\frac{x^6}{6}\\
&=&1+\frac{1}{1!}x^2+\frac{1}{2!}\left(x^2\right)^2+\frac{1}{3!}\left(x^2\right)^3\\
&=&\displaystyle\sum_{i=0}^{3} \frac{1}{i!}\left(x^2\right)^i.
\end{eqnarray*}
Consequently, it can be inferred that
\begin{equation}
y_n(x)=\displaystyle\sum_{i=0}^{n} \frac{1}{i!}\left(x^2\right)^i, n\in \N,
\end{equation}
which tends to $e^{x^2}$ as $n\rightarrow \infty.$
\par  The exact solution is known to be $y(x)=e^{x^2}.$ Figure \ref{sumud33} displays the exact and approximate solutions for Example \ref{sumud11}.
\begin{figure}
\includegraphics[width=12.0cm ,height=10.0cm]{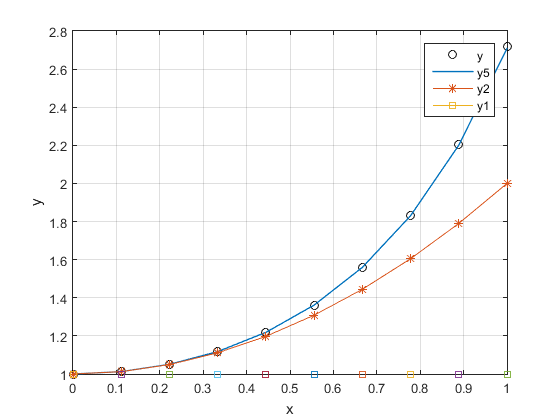}
 \caption{Graph of $y(x)=e^{x^2}.$}
 \label{sumud33}
\end{figure}
\end{example}
\subsection{Volterra integro-differential equation}
Physics, engineering and biology books illustrate the application of the Volterra integro-differential equations (VIDEs) in mathematical modeling of various processes such as diffusion, glassforming, heat transfer, biological species coexisting and nanohydrodynamics. The approximate or numerical solutions of VIDEs have been obtained by using some methods which include adomian decomposition method \cite{Momani}, collocation method \cite{Rawashdeh, Ma1}, fractional differential transform method \cite{Arikoglu}, variational iteration method \cite{Nawaz, Odibat, Sayevand}, reproducing kernel method \cite{Wei} and wavelet method \cite{Saeedi2, Saeedi1, Zhu2, Zhu1, Wang2, Wang1}. The efficiency of the STM method for obtaining the solutions of VIDEs is presented in this paper. An example is given to illustrate the procedures. Each succeeding term after the first iteration is derived by truncating the HOT while computing each subsequent iteration. The graph is shown by using Matlab to compare the exact solution with the solution obtained via the STM.
 \begin{example}\label{sumud40}
Consider first order nonlinear VIDE,
 \begin{eqnarray}\label{sumud34}
 y'(x)-\int^{x}_{0}y^2\left(\frac{t}{2}\right)dt=0,
 \end{eqnarray}
 for $x\geq 0$ and which is subject to the initial condition $y(0) = 1.$
\par {\bf Solution:} Taking the derivative of (\ref{sumud34}) gives
 \begin{eqnarray}\label{sumud35}
 y''(x)-y^2\left(\frac{x}{2}\right)=0.
 \end{eqnarray}
 Substituting $x=0$ in VIDE (\ref{sumud34}) gives $y'(0)=1.$
 Then the VIDE (\ref{sumud34}) subject to initial condition $y(0) = 1$ is reduced to pantograph-type equation (\ref{sumud35})
subject to initial conditions
 \begin{eqnarray}\label{sumud36}
 y(0)&=&1 \nonumber\\
  y'(0)&=&1.
 \end{eqnarray}
 Taking the ST of (\ref{sumud35}) gives
 \begin{equation}
\frac{Y(u)}{u^2}-\frac{y(0)}{u^2}-\frac{y'(0)}{u}=S\left[y^2\left(\frac{x}{2}\right)\right].
\end{equation}
Since $y(0)=y'(0)=1,$ it yields that
\begin{equation}
\frac{Y(u)}{u^2}-\frac{1}{u^2}-\frac{1}{u}=S\left[y^2\left(\frac{x}{2}\right)\right].
\end{equation}
Then the variational iteration formula emerges as
\begin{eqnarray}\label{sumud37}
Y_{n+1}(u)&=&Y_n(u)+{\varphi}(u)\left(\frac{Y_n(u)}{u^2}-\frac{1}{u^2}-\frac{1}{u}-S\left[y_n^2\left(\frac{x}{2}\right)\right]\right).
\end{eqnarray}
 $S\left[y_n^2\left(\frac{x}{2}\right)\right]$ is the restricted variation and taking the classical variation operator on both sides of (\ref{sumud37}) gives the Lagrange multiplier as
\begin{equation}
{\varphi}(u)=-u^2.
\end{equation}
Substitute for ${\varphi}(u)$ in (\ref{sumud37}) and take its inverse-Sumudu transform $S^{-1}.$ It brings forth the explicit iteration formula 
\begin{eqnarray}\label{sumud38}
y_{n+1}(x)&=&y_n(x)+S^{-1}\left[-u^2\left(\frac{Y_n(u)}{u^2}-\frac{1}{u^2}-\frac{1}{u}-S\left[y_n^2\left(\frac{x}{2}\right)\right]\right)\right]\nonumber\\
&=&y_1(x)+S^{-1}\left[u^2\left(S\left[y_n^2\left(\frac{x}{2}\right)\right]\right)\right],
\end{eqnarray}
with the initial approximation which is given by
$$y_1(x)=S^{-1}\left[-u^2(\frac{-1}{u^2}-\frac{-1}{u})\right]=S^{-1}[1+u]=1+x.$$
\begin{eqnarray*}
y_{2}(x)&=&y_1(x)+S^{-1}\left[u^2\left(S\left[y_1^2\left(\frac{x}{2}\right)\right]\right)\right]\\
&=&1+x+S^{-1}\left[u^2\left(S\left[\left(1+\frac{x}{2}\right)^2\right]\right)\right]\\
&=&1+x+S^{-1}\left[u^2\left(S\left[1\right]\right)\right]~(\mbox{Truncate HOT})\\
&=&1+x+S^{-1}\left[u^2\right]\\
&=&1+x+\frac{x^2}{2!}.
\end{eqnarray*}
Next is to evaluate $y_3(x)$ to obtain the next term.
\begin{eqnarray*}
y_{3}(x)&=&y_1(x)+S^{-1}\left[u^2\left(S\left[y_2^2\left(\frac{x}{2}\right)\right]\right)\right]\\
&=&1+x+S^{-1}\left[u^2\left(S\left[\left(1+\frac{x}{2}+\frac{x^2}{8}\right)^2\right]\right)\right]\\
&=&1+x+S^{-1}\left[u^2\left(S\left[1+x\right]\right)\right]~(\mbox{Truncate HOT})\\
&=&1+x+S^{-1}\left[u^2\left(1+u\right)\right]\\
&=&1+x+S^{-1}\left[u^2+u^3\right]\\
&=&1+x+\frac{x^2}{2!}+\frac{x^3}{3!}.
\end{eqnarray*}
Next is to evaluate $y_4(x)$ to obtain the next term.
\begin{eqnarray*}
y_{4}(x)&=&y_1(x)+S^{-1}\left[u^2\left(S\left[y_2^2\left(\frac{x}{2}\right)\right]\right)\right]\\
&=&1+S^{-1}\left[u^2\left(S\left[\left(1+\frac{x}{2}+\frac{x^2}{8}+\frac{x^3}{48}\right)^2\right]\right)\right]\\
&=&1+x+S^{-1}\left[u^2\left(S\left[\left(1+x+\frac{x^2}{2}\right)\right]\right)\right]~(\mbox{Truncate HOT})\\
&=&1+x+S^{-1}\left[u^2\left(1+u+u^2\right)\right]\\
&=&1+x+S^{-1}\left[u^2+u^3+u^4\right]\\
&=&1+x+\frac{x^2}{2!}+\frac{x^3}{3!}+\frac{x^4}{4!}\\
&=&\displaystyle\sum_{i=0}^{4} \frac{1}{i!}x^i.
\end{eqnarray*}
Accordingly, it can be deduced that
\begin{equation}
y_n(x)=\displaystyle\sum_{i=0}^{n} \frac{1}{i!}x^i, n\in \N,
\end{equation}
which inclines to $e^x$ as $n\rightarrow \infty.$
\par  The exact solution is known to be $y(x)=e^x.$ Figure \ref{sumud39} displays the exact and approximate solutions for Example \ref{sumud40}.
\begin{figure}
\includegraphics[width=12.0cm ,height=10.0cm]{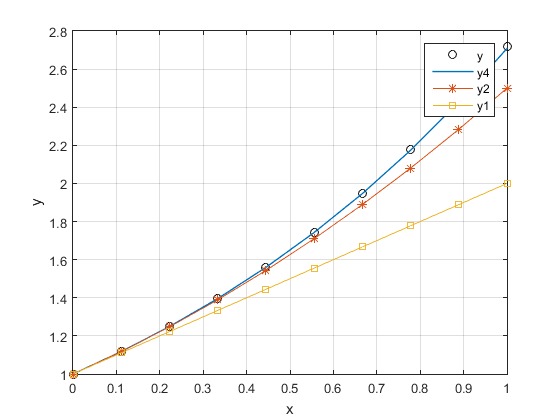}
 \caption{Graph of $y(x)=e^x.$}
 \label{sumud39}
\end{figure}
\end{example}
\vspace{1.0cm}
\noindent{\bf Conclusion}: This paper shows how the Sumudu Transform Method is used to identify the Lagrange multipliers for nonlinear DDEs with variable coefficients. It displays the craft of Sumudu Transform Method in reducing the complex computational work as compared to the popular methods. The Sumudu Transform Method is used to find the exact and approximate solutions of DDEs of pantograph type with variable coefficients and nonlinear Volterra integro-differential equations of pantograph type. This study generalises and extends the previous works in the literature where pantograph type equations with constant coefficients were considered (See e.g, \cite{Vilu1}). The graphs of the expository examples which were given, show that the solutions obtained via the STM are in good agreement with the exact solutions.\\
 
\noindent {\bf List of Abbreviations}\\
ST: Sumudu Transform\\
STM: Sumudu Transform Method \\
DDEs: Delay Differential Equations\\
VIM: Variational iteration method\\
PTEs: Pantograph Type Equations\\ 
HOT: Higher Other Terms\\
VIDEs: Volterra integro-differential equations\\ 

\noindent {\bf Acknowledgements}: The first author acknowledges with thanks the postdoctoral fellowship and financial support from the DSI-NRF Center of Excellence in Mathematical and Statistical Sciences (CoE-MaSS). Opinions expressed and conclusions arrived are those of the authors and are not necessarily to be attributed to the CoE-MaSS.

\footnotesize
\noindent {\bf Conflicts of Interest}\\
The authors declare no conflict of interest.


\end{document}